\documentclass[12pt,a4paper]{article}
\usepackage[ansinew]{inputenc}
\usepackage{latexsym,amsfonts,amssymb,theorem,amsmath}

\usepackage{url}
\newenvironment{proof}{\noindent{\bf Proof.}}{\ignorespaces\rule{1pt}{0pt}\hfill$\Box$\medskip}
\newtheorem{theorem}{Theorem}[section]
\newtheorem{lemma}[theorem]{Lemma}
\newtheorem{cor}[theorem]{Corollary}

\newtheorem{remark}[theorem]{Remark}
{\theorembodyfont{\rm}}
{\theorembodyfont{\rm}}

\newcommand{\PG}{{\rm PG}}

\newcommand{\pg}[1]{\ensuremath{\mathrm{PG}(#1,q)}}
 
\newcommand{\Q}[1]{\ensuremath{\mathrm{Q}(#1,q)}}                 
\newcommand{\Qplus}[1]{\ensuremath{\mathrm{Q^+}(#1,q)}}           
\newcommand{\Qmin}[1]{\ensuremath{\mathrm{Q^-}(#1,q)}}            

\title{On the small weight codewords of the functional codes $C_2(\mathrm{Q})$, $\mathrm{Q}$ a non-singular  quadric}
\author{F. A. B. Edoukou\and A. Hallez\thanks{The research of this author is supported by a research grant of the Research council of Ghent University.} \and F. Rodier\and L. Storme}

\begin{document}

\maketitle

\begin{abstract} 
We study the small weight codewords of the functional code $C_2(\mathrm{Q})$, with $\mathrm{Q}$ a non-singular quadric of $\pg{N}$. We prove that the small weight codewords correspond to the intersections of $\mathrm{Q}$ with the singular quadrics of $\pg{N}$ consisting of two hyperplanes. We also calculate the number of codewords having these small weights.
\end{abstract}
\section{Introduction}

Consider a  non-singular quadric $\mathrm{Q}$ of $\pg{N}$. Let $\mathrm{Q}=\{P_1,\ldots,P_n\}$, where we normalize the coordinates of the points $P_i$ with respect to the leftmost non-zero coordinate. Let ${\mathcal F}$ be the set of all homogeneous quadratic polynomials $f(X_0,\ldots,X_N)$ defined by $N+1$ variables. The functional code 
$C_2(\mathrm{Q})$ is the linear code 
\[C_2(\mathrm{Q})=\{(f(P_1),\ldots,f(P_n))|| f\in \mathcal{F}\cup\{0\}\},\]
defined over $\mathbb{F}_q$. 

This linear code has length $n=|\mathrm{Q}|$ and  dimension $k=\left(\begin{array}{c}N+2\\2\end{array}\right)-1$. The third fundamental parameter of this linear code $C_2(\mathrm{Q})$ is its minimum distance $d$. 

We determine the 5 or 6 smallest weights of $C_2(\mathrm{Q})$ via geometrical arguments. Every homogeneous quadratic polynomial $f$ in $N+1$ variables defines a quadric $\mathrm{Q}':f(X_0,\ldots,X_N)=0$. The small weight codewords of $C_2(\mathrm{Q})$ correspond to the quadrics of $\pg{N}$ having the largest intersections with $\mathrm{Q}$. 

We prove that these small weight codewords correspond to quadrics $\mathrm{Q}'$ which are the union of two hyperplanes of $PG(N,q)$. Since there are different possibilities for the intersection of two hyperplanes with a non-singular quadric, we determine in this way the 5 or 6 smallest weights of the functional code $C_2(\mathrm{Q})$.

We also determine the exact number of codewords having the 5 or 6 smallest weights.
 \section{Quadrics in $\pg{N}$}
 
 The non-singular quadrics in $\pg{N}$ are equal to:
 
 \begin{itemize}
 \item the non-singular parabolic quadrics $\mathrm{Q}(2N,q)$ in $\pg{2N}$ having standard equation $X_0^2+X_1X_2+\cdots+X_{2N-1}X_{2N}=0$. These quadrics contain $q^{2N-1}+\cdots+q+1$ points, and the largest dimensional spaces contained in a non-singular parabolic quadric of $\pg{2N}$ have dimension $N-1$,
 \item the non-singular hyperbolic quadrics $\mathrm{Q}^+(2N+1,q)$ in $\pg{2N+1}$ having standard equation $X_0X_1+\cdots+X_{2N}X_{2N+1}=0$. These quadrics contain $(q^N+1)(q^{N+1}-1)/(q-1)=q^{2N}+q^{2N-1}+\cdots+q^{N+1}+2q^N+q^{N-1}+\cdots+q+1$ points, and the largest dimensional spaces contained in a non-singular hyperbolic quadric of $\pg{2N+1}$ have dimension $N$,
  \item the non-singular elliptic quadrics $\mathrm{Q}^-(2N+1,q)$ in $\pg{2N+1}$ having standard equation $f(X_0,X_1)+X_2X_3+\cdots+X_{2N}X_{2N+1}=0$, where $f(X_0,X_1)$ is an irreducible quadratic polynomial over $\mathbb{F}_q$. These quadrics contain $(q^{N+1}+1)(q^N-1)/(q-1)=q^{2N}+q^{2N-1}+\cdots+q^{N+1}+q^{N-1}+\cdots+q+1$ points, and the largest dimensional spaces contained in a non-singular elliptic quadric of $\pg{2N+1}$ have dimension $N-1$.
  \end{itemize}
  
  All the quadrics of $\pg{N}$, including the non-singular quadrics, can be described as a quadric having an $s$-dimensional vertex $\pi_s$ of singular points, $s\geq -1$, and having a non-singular base $\mathrm{Q}_{N-s-1}$ in an $(N-s-1)$-dimensional space  skew to $\pi_s$, denoted by $\pi_sQ_{N-s-1}$.
  
  We denote the largest dimensional spaces contained in a quadric by the {\em generators} of this quadric.
  
  Since we will make heavily use of the sizes of (non-)singular quadrics of $\pg{N}$, we list these sizes explicitly.
  
  \begin{itemize}
  \item In $\pg{N}$, a quadric having an $(N-2d-2)$-dimensional vertex and a hyperbolic quadric $\mathrm{Q}^+(2d+1,q)$ as base has size 
  \[q^{N-1}+\cdots+q^{N-d}+2q^{N-d-1}+q^{N-d-2}+\cdots+q+1.\]
  \item In $\pg{N}$, a quadric having an $(N-2d-2)$-dimensional vertex and an elliptic quadric $\mathrm{Q}^-(2d+1,q)$ as base has size 
  \[q^{N-1}+\cdots+q^{N-d}+q^{N-d-2}+\cdots+q+1.\]
  \item In $\pg{N}$, a quadric having an $(N-2d-1)$-dimensional vertex and a parabolic quadric $\mathrm{Q}(2d,q)$ as base has size 
  \[q^{N-1}+q^{N-2}+\cdots+q+1.\]
  \end{itemize}
  
  We note that the size of a (non-)singular quadric having a non-singular hyperbolic quadric as base, is always larger than the size of a (non-)singular quadric having a non-singular parabolic quadric as base, which is itself always larger than the size of a (non-)singular quadric having a non-singular elliptic quadric as base.
  
  The quadrics having the largest size are the union of two distinct hyperplanes of $\pg{N}$, and have size $2q^{N-1}+q^{N-2}+\cdots+q+1$. The second largest quadrics in $\pg{N}$ are the quadrics having an $(N-4)$-dimensional vertex and a non-singular 3-dimensional hyperbolic quadric $\mathrm{Q}^+(3,q)$ as base. These quadrics have size $q^{N-1}+2q^{N-2}+q^{N-3}+\cdots+q+1$. The third largest quadrics in $\pg{N}$ have an $(N-6)$-dimensional vertex and a non-singular hyperbolic quadric $\mathrm{Q}^+(5,q)$ as base. These quadrics have size $q^{N-1}+q^{N-2}+2q^{N-3}+q^{N-4}+\cdots+q+1$.
  
  As we mentioned in the introduction, the smallest weight codewords of the code $C_2(\mathrm{Q})$ correspond to the  largest intersections of $\mathrm{Q}$ with other quadrics $\mathrm{Q}'$ of $\pg{N}$. Let $V$ be the intersection of the  quadric $\mathrm{Q}$ with the quadric $\mathrm{Q}'$. Two distinct quadrics $\mathrm{Q}$ and $\mathrm{Q}'$ define a unique pencil of quadrics $\lambda \mathrm{Q}+\mu \mathrm{Q}'$, $(\lambda,\mu)\in \mathbb{F}_q^2\setminus \{(0,0)\}$.

Let $V=\mathrm{Q}\cap \mathrm{Q}'$, then $V$ also lies in every quadric $\lambda \mathrm{Q}+\mu \mathrm{Q}'$ of the pencil of quadrics defined by $\mathrm{Q}$ and $\mathrm{Q}'$. A large intersection implies that there is a large quadric in the pencil. The $q+1$ quadrics of the pencil contain altogether $\vert \pg{N}\vert+q\vert V\vert$ points, since the points of $V$ lie in all the quadrics of the pencil and the other points of $\pg{N}$ lie in exactly one such  quadric. So there is a quadric in the pencil containing at least $(\vert \pg{N}\vert+q\vert V\vert)/(q+1)$ points. 

If there is a quadric in the pencil which is equal to the union of two hyperplanes, then we are at the desired conclusion that the largest intersections of $\mathrm{Q}$ arise from the intersections of $\mathrm{Q}$ with the quadrics which are the union of two hyperplanes. So assume that all $q+1$ quadrics in this pencil defined by $\mathrm{Q}$ and $\mathrm{Q}'$ are irreducible; we try to find a contradiction. As already mentioned above, the largest irreducible quadrics are cones with vertex \pg{N-4} and base \Qplus{3}, and the second largest irreducible quadrics are cones with vertex \pg{N-6} and base \Qplus{5}.

\begin{theorem}\label{thm:0}
In $\pg{N}$, with $N \geqslant 6$, or $N=5$ and $\mathrm{Q}=\Qmin{5}$, if $\vert V \vert > q^{N-2}+3q^{N-3}+3q^{N-4}+2q^{N-5}+\cdots+2q+1$, then in the pencil of quadrics defined by $\mathrm{Q}$ and $\mathrm{Q}'$, there is a quadric consisting of two hyperplanes. 
\end{theorem}
\begin{proof}
Suppose that there is no quadric consisting of two hyperplanes in the pencil of quadrics.

If $\vert V \vert > q^{N-2}+2q^{N-3}+2q^{N-4}+q^{N-5}+\cdots+q+1$, then $(\vert \pg{N}\vert+q\vert V\vert)/(q+1)>| \pi_{N-6} \Qplus{5}|$, so there is a singular quadric  $\pi_{N-4} \Qplus{3}$ in the pencil of quadrics. 

With the lines of one regulus of $\Qplus{3}$, together with $\pi_{N-4}$, we form $q+1$ different $(N-2)$-spaces $\pi_{N-2}$. We wish to have that at least one of these $(N-2)$-spaces intersects $\mathrm{Q}$ in two $(N-3)$-dimensional spaces. All points of $V$ appear in at least one of these $(N-2)$-dimensional spaces $\Pi_{N-2}$, so for some space $\pi_{N-2}$, we have that $\vert\pi_{N-2}\cap V\vert\geqslant \vert V\vert/(q+1)$.
 
If $\vert V\vert/(q+1)> \vert  \pi_{N-6}\Qplus{3}\vert$, then $\pi_{N-2}\cap \mathrm{Q}$ is the union of two $(N-3)$-spaces.
When $\vert V\vert>q^{N-2}+3q^{N-3}+3q^{N-4}+2q^{N-5}+\cdots+2q+1$, then this is valid.
So $\pi_{N-2}\cap \mathrm{Q}=\pi^1_{N-3}\cup\pi^2_{N-3}$.

These two $(N-3)$-dimensional spaces are contained in $V$, so belong to $\mathrm{Q}$. 
This means that $\mathrm{Q}$ must have subspaces of dimension $N-3.$ The next table shows that this can only occur in small dimensions.

\begin{center}
\begin{tabular}{|c|c|c|}\hline
quadric&dimension generator& property fulfilled\\
\hline
$\mathrm{Q}$=\Qplus{N=2N'+1}&$N'$&$N'\leq 2$\\ \hline
$\mathrm{Q}$=\Qmin{N=2N'+1}&$N'-1$&$N'\leq 1$\\ \hline
$\mathrm{Q}$=\Q{N=2N'}&$N'-1$&$N'\leq 2$ \\ \hline
\end{tabular}
\end{center}
Except for the small cases for $N'$, we have a contradiction, so there is a quadric consisting of two hyperplanes in the pencil of quadrics defined by $\mathrm{Q}$ and $\mathrm{Q}'$.
\end{proof}

\begin{remark}
First of all we say something about the sharpness of the bound in Theorem \ref{thm:0}. Therefore we refer to \cite{AC:95}.
In a pencil of $q+1$ non-singular elliptic quadrics $\mathrm{Q}^-(N,q)$ not containing hyperplanes, the size of the intersection of 2 quadrics is:

\begin{eqnarray*}
\vert Q_1\cap Q_2\vert&=&q^{N-2}+q^{N-3}+\cdots+q^\frac{N+1}{2}+q^\frac{N-5}{2}+\cdots+q+1
\end{eqnarray*}
\end{remark}

Since the problem is now solved for dimensions $N$ up to 4 \cite{FE:07,FE:08}, there is only one case still open.
From now on, $\mathrm{Q}$ will be the hyperbolic quadric \Qplus{5}.

If $\vert V \vert > q^3+2q^2+2q+1$, then there is a singular quadric $ \pi_{N-4}\Qplus{3}= L\Qplus{3}$ in the pencil of quadrics, if we assume that there is no quadric in the pencil which is the union of two hyperplanes.

We form solids $\omega_1,\ldots,\omega_{q+1}$ with $L$ and the lines of one regulus of the base \Qplus{3}. If $\vert V\vert >q^3+3q^2+3q+1$, $\vert V\vert/(q+1)>\vert  \pi_{N-6}\Qplus{3}\vert$,  so there is a solid through $L$ of $L\Qplus{3}$ intersecting $\mathrm{Q}$ in two planes.
 
Now we have three different cases:

\begin{enumerate}
\item $L\subset V$,
\item $\vert L \cap V \vert =1$,
\item $\vert L \cap V \vert =2$.
\end{enumerate}

\begin{lemma}\label{lem:1}
For \Qplus{5}, if $\vert V \vert >q^3+4q^2+1$ and $L\subset V$, then there is a quadric consisting of two hyperplanes in the pencil of quadrics defined by $\mathrm{Q}$ and $\mathrm{Q}'$.
\end{lemma}
\begin{proof}
Assume that  no quadric in the pencil is the union of two hyperplanes. Then we have already a singular quadric $ L\Qplus{3}$ in the pencil and there is a solid $\omega_1$ through $L$ intersecting $\mathrm{Q}$ in 2 planes. Now $L$ lies in one or both of these  planes, since $L \subset V$.

Every point of $V$ lies in at least one of the $q+1$ solids $\omega_1,\ldots,\omega_{q+1}$ through $L$. Now
$$\vert V \vert-\mbox{(union of 2 planes)}>q^3+4q^2+1-(2q^2+q+1)=q^3+2q^2-q.$$

So one of the $q$ remaining solids of $ \omega_2,\ldots,\omega_{q+1}$ contains at least $$\vert L\vert+\frac{q^3+2q^2-q}{q}=q^2+3q$$points.

So one solid $\omega_2$ contains more  than $|\Qplus{3}|$ points of $V$, so $\omega_2$ intersects $\mathrm{Q}$ in the  union of two planes. One of these planes contains $L$, so $L$ lies already in two planes of \Qplus{5}.

Now one of the $q-1$ remaining solids $\omega_3,\ldots,\omega_{q+1}$ contains more than $$q+1+(q^3+2q^2-q-2q^2)/(q-1)=q^2+2q+1$$ points of $V$.

Again this implies that there is a solid $\omega_3$ intersecting $\mathrm{Q}$  in  the union of two planes, with at least one of them containing $L$.
This gives us at least three planes of \Qplus{5} through $L$, which is impossible. We have a contradiction.
So there is a quadric consisting of 2 hyperplanes in the pencil of quadrics defined by $\mathrm{Q}$ and $\mathrm{Q}'$.  
\end{proof}

\begin{lemma}
For \Qplus{5}, if $\vert V \vert >q^3+5q^2+1$, then the case $|L\cap V|=1$ does not occur.
\end{lemma}
\begin{proof}
Assume that  no quadric in the pencil of $\mathrm{Q}$ and $\mathrm{Q}'$ is the union of two hyperplanes. Then we have already a singular quadric $ L\Qplus{3}$ in the pencil of quadrics. In this quadric,  the line $L$ is skew to the solid of \Qplus{3}.

But $L$ is a tangent line to \Qplus{5} in a point $R$ since $L$ is contained in the cone $L$\Qplus{3}, but $L$ shares only one point with \Qplus{5}.

Using the same arguments as in the preceding lemma, we prove that at least three solids defined by the line $L$ and lines of one regulus of the base \Qplus{3} intersect $\mathrm{Q}$ in two planes. These planes all pass through $R$, so they lie in the tangent hyperplane $T_R(\mathrm{Q})$, which intersects $\mathrm{Q}$ in a cone with vertex $R$ and base $\Qplus{3}'$.  Two such planes of $V$ in the same solid of $L$\Qplus{3} through $L$ intersect in a line, so they define lines of the opposite reguli of the base $\Qplus{3}'$ of this tangent cone. This shows that the 4-space defined by $R$ and the base \Qplus{3}'
shares already six planes with $\mathrm{Q}$. By B\'ezout's theorem, the cone $R$\Qplus{3}' is contained in $V$.

Consider a hyperplane through $L$; this intersects $L\Qplus{3}$ either in a cone $L\mathrm{Q}(2,q)$ or in the union of two solids. So the tangent hyperplane $T_R(\mathrm{Q})$ cannot intersect $L\Qplus{3}$ in a cone $R\Qplus{3}'$.

This gives us a contradiction.
\end{proof}

\begin{lemma}
For \Qplus{5}, if $\vert V \vert >q^3+5q^2-q+1$ and $|L\cap V|=2$, then there is a quadric consisting of two hyperplanes in the pencil of quadrics defined by $\mathrm{Q}$ and $\mathrm{Q}'$.
\end{lemma}
\begin{proof}
Assume that  no quadric in the pencil defined by $\mathrm{Q}$ and $\mathrm{Q}'$ is the union of two hyperplanes. Then we have already a singular quadric $ L\Qplus{3}$ in the pencil and there is a solid $\omega_1$ through $L$ intersecting $\mathrm{Q}=\Qplus{5}$ in two planes. Assume that $L\cap V=\{R, R'\}$. Let $\mathrm{Q}^+(3,q)_L$ be the polar quadric of $L$ w.r.t. $\mathrm{Q}^+(5,q)$ and let $\mathrm{Q}^+(3,q)_L$ lie in the solid $\pi_3$.

By the same counting arguments as in Lemma \ref{lem:1}, we know that if $\vert V \vert >q^3+5q^2-2q+2$, then there are 3 solids $\langle L,L_i\rangle$, with $i=1,2,3$, and all $L_i$ belonging to the same regulus of \Qplus{3}, intersecting $\mathrm{Q}$ in 2 planes.  
For every solid $\langle L,L_i\rangle$, we denote by $\tilde{L_i}$ the line that the 2 planes have in common, and $\pi_{i1}=\langle R, \tilde{L_i}  \rangle , \pi_{i2}=\langle R', \tilde{L_i}  \rangle$.  Then
$\tilde{L_i}=\pi_{i1}\cap\pi_{i2}\subset R^\bot\cap R'^\bot=\pi_3$, with $\bot$  the polarity w.r.t. $\Qplus{5}$.
We use the same arguments for the opposite regulus.  This gives us again 3 solids $\langle L,M_i\rangle$, $i=1,2,3$, intersecting $\mathrm{Q}$ in 2 planes. We denote by $\tilde{M_i}$ the line in the intersection of these 2 planes.

These lines $\tilde{L_i}$ and $\tilde{M_i}$ belong to the hyperbolic quadric $\Qplus{3}_L$ in $R^\bot\cap R'^\bot$, which is the basis for $ R\Qplus{3}_L $ as well as for $ R'\Qplus{3}_L $.
The quadric $ R \Qplus{3}_L$ shares 6 planes with $ L \Qplus{3}$. By B\'ezout's theorem,  if $ R \Qplus{3}_L \not \subset  L \Qplus{3} $, then the intersection would be of degree 4, so $ R \Qplus{3}_L \subset  L \Qplus{3} \cap \mathrm{Q}$. Similarly, $ R' \Qplus{3}_L \subset  L \Qplus{3} \cap \mathrm{Q}$.

The cone $L \Qplus{3}$ intersects $\mathrm{Q}$ in 2 tangent cones $ R \Qplus{3}_L$ and $ R' \Qplus{3}_L$. 
We will now look at the pencil of quadrics defined  by $\mathrm{Q}$ and $ L \Qplus{3}=\mathrm{Q}'$. 

Let $P$ be a point of $\pi_3\backslash\Qplus{3}_L$. The points of $\pg{5}\backslash(\mathrm{Q}\cap \mathrm{Q}')$ lie in exactly one quadric of the pencil defined by $\mathrm{Q}$ and $\mathrm{Q}'$. For the point $P$, this must be the quadric consisting of the two hyperplanes $ \langle R,\pi_3\rangle$ and $ \langle R',\pi_3  \rangle$. For  $ \langle R, \pi_3\rangle$ contains a cone $R  \Qplus{3}_L$ and the point $P$ of this quadric, so this is one point too much for a quadric. 

So one quadric of the pencil consists of the union of 2 hyperplanes. 
\end{proof}

\begin{cor}
For $\Qplus{5}$, if $\vert V \vert >q^3+5q^2+1$, then the intersection  of $\Qplus{5}$ with the other quadric $\mathrm{Q}'$ is equal to the intersection of $\Qplus{5}$ with the union of two hyperplanes.
\end{cor}
\section{Dimension 4}

We consider a pencil of quadrics $\lambda\mathrm{Q}+\mu\mathrm{Q}'$ in PG$(4,q)$, with Q a non-singular parabolic quadric $\mathrm{Q}(4,q)$. Let $V=\mathrm{Q}\cap \mathrm{Q}'$. If $\vert V\vert>q^2+q+1$, then there is at least one cone $P\Qplus{3}$ in this pencil. 

\begin{lemma}\label{lem:2}
If $\vert V\vert >q^2+(x+1)q+1$, then $x$ planes through $P$ of the same regulus of $P\Qplus{3}$ intersect Q in 2 lines.
\end{lemma}
\begin{proof}
Consider one regulus of  $P\Qplus{3}$. We wish to have that $x$ planes $PL$, with $L$ a line of this regulus, intersect Q in 2 lines. So for the first plane,  this means that $\frac{\vert V\vert}{q+1}>q+1$, since every point of $V$ lies in one of the $q+1$ planes $PL$. 
For the $x$-th plane, we have already $x-1$ planes which intersect Q in 2 lines. We impose that $\frac{\vert V\vert-(x-1)(2q+1)}{q-x+2}>q+1$ to guarantee that the $x$-th plane also intersects Q in 2 lines. This reduces to  $\vert V\vert >q^2+(x+1)q+1$.  
\end{proof}

Denote by $L_i$ the lines of one regulus of \Qplus{3} and by $M_i$ the lines of the opposite regulus of \Qplus{3}, with $i=1,2,\ldots,q+1$. Denote by $l_{i1}, l_{i2}$, resp. $m_{i1}, m_{i2}$, the lines of $\mathrm{Q}\cap PL_i$, resp. $\mathrm{Q}\cap PM_i$.\\

We have to look at 2 cases now, whether $P\in V$ or whether   $P\not\in V$.\\

\textbf{CASE I: $P\in V$}

\begin{theorem}\label{thm:1}
For \Q{4}, if $\vert V\vert>q^2+6q+1$ and $P\in V$, then V consists of the union of a cone $P\Q{2}$ and another 3-dimensional quadric.
\end{theorem}
\begin{proof}
If we consider one regulus of the base $P\Qplus{3}$, then, by the preceding lemma, there are $x\geqslant5$ planes each containing 2 lines of $V$, of   which at least one goes through $P$. This gives us at least $ x\geqslant 5$ lines through $P$ in $\mathrm{Q}(4,q)\cap P\Qplus{3}$. These $x$ lines lie on the tangent cone $P\mathrm{Q}(2,q)$ in  $T_P(\mathrm{Q}(4,q))$. By B\'ezout's theorem, since $x\geqslant5$, this cone $P\mathrm{Q}(2,q)$ lies completely in $\mathrm{Q}(4,q)$ and in $P\Qplus{3}$. \\
Since $\mathrm{Q}(4,q)\cap P\Qplus{3}$ is an algebraic variety of degree 4 and dimension 2, and since  $\vert V\vert > \vert P\mathrm{Q}(2,q)\vert$, $V$ is the union of $P\mathrm{Q}(2,q)$ and another 3-dimensional quadric.
\end{proof}

\textbf{CASE II: $P\not\in V$}

\begin{theorem}\label{thm:2}
For \Q{4}, if $\vert V\vert>q^2+11q+1$ and $P\not\in V$, then for $q>7$, $V$ consists of the union of 2 hyperbolic quadrics. 
\end{theorem}
\begin{proof}
We use the notations introduced after the proof of Lemma \ref{lem:2}.

Without loss of generality,  we can assume that the lines of $PL_i$ lying on Q intersected by $m_{11}$ (resp. $m_{12}$) are the lines $l_{i1}$ (resp. $l_{i2}$), $i=1,\ldots,x$. So $m_{11}$ and $m_{12}$ are both intersected by $x$ lines of Q. 
 
The line $m_{21}$ will intersect at least $\lceil \frac{x}{2}\rceil$ of the lines $l_{i1}$. This means that $m_{21}$ has these transversals in common with $m_{11}$. Assume that these lines are the lines $l_{11}, \cdots, l_{\lceil\frac{x}{2}\rceil 1}$.
Also $m_{31}$ has at least $\lceil \frac{x}{2}\rceil$ transversals in common with $m_{11}$. 

Assume that at least 2 of those transversals also intersect $m_{21}$, then $m_{11},m_{21}, m_{31}$ define a 3-dimensional hyperbolic quadric \Qplus{3} sharing 5 lines with Q$(4,q)$.

Otherwise, at least $x-1$ transversals out of the $x$ selected transversals to $m_{11}$ are intersecting one of $m_{21}$ and $m_{31}$. 
Suppose now that $m_{41}$ shares  at least $ \lceil \frac{x}{2}\rceil$ transversals with $m_{11}$. One of them could be skew to $m_{21}$ and $m_{31}$, but at least $\lceil\frac{x}{2}\rceil-1$ of them intersect $m_{21}$ or $m_{31}$. At least $\frac{\frac{x}{2}-1}{2}$ of them intersect, for instance, $m_{21}$. If this is at least 2, then  $m_{11},m_{21}, m_{41}$ define a 3-dimensional hyperbolic quadric \Qplus{3} sharing $5$ lines with $\mathrm{Q}(4,q)$. Therefore, we obtain the same conclusion that $V$ contains a 3-dimensional hyperbolic quadric when $x\geqslant 10$. Lemma \ref{lem:2} implies that we need to impose that $\vert V\vert > q^2+11q+1$.
Since in both cases, there is a 3-dimensional  hyperbolic quadric \Qplus{3} sharing $5$ lines with $\mathrm{Q}(4,q)$,   B\'ezout's theorem implies  that $\Qplus{3}\subset \mathrm{Q}(4,q)$. So $V$ consists of  \Qplus{3} and another 3-dimensional quadric. The remaining lines of $V$ are $10$ skew lines of planes $PL_i$ and $10$ skew lines of planes $PM_j$, and these  lines of $V$ lying in  $PL_i$ intersect the lines of $V$ lying in  $PM_j$. So these lines also form a 3-dimensional hyperbolic quadric \Qplus{3}.  
\end{proof}


\begin{theorem}
For \Q{4}, if $\vert V\vert>q^2+11q+1$, then there is a union of 2 hyperplanes in the pencil of quadrics defined by Q and Q'.
\end{theorem}
\begin{proof}
By Theorems \ref{thm:1} and \ref{thm:2},  $V$ consists of a 3-dimensional hyperbolic quadric $\Qplus{3}$ in a solid $\pi_3$ and another 3-dimensional quadric.
Let $R$ be a point of $\pi_3\backslash V$. The points of $\pg{4}\backslash(\mathrm{Q}\cap \mathrm{Q}')$ lie in exactly one quadric of the pencil. Let Q$''$ be the unique quadric in the pencil defined by Q and Q$'$ containing $R$. So $ \pi_3$ shares with Q$''$ a quadric and an extra point $R$, so this is one point too much for a quadric, hence there is a quadric in the pencil defined by $\mathrm{Q}$ and $\mathrm{Q}'$ containing a hyperplane, so a quadric in the pencil defined by two hyperplanes. 
\end{proof}
\section{Tables}

For the standard properties and notations on quadrics, we refer to \cite{JWPH:91}.

\subsection{The hyperbolic quadric in PG$(2l+1,q)$}
We know that the largest intersections of a non-singular hyperbolic quadric $\Qplus{2l+1}$ in $\PG(2l+1,q)$ with the other quadrics are the intersections of $\Qplus{2l+1}$ with the quadrics which are the union of two hyperplanes $\Pi_1$ and $\Pi_2$. We now discuss all the different possibilities for the intersections of $\Qplus{2l+1}$ with the union of two hyperplanes. This then gives the five or six smallest weights of the functional codes $C_2(\Qplus{2l+1})$, and the numbers of  codewords having these weights.

We  start the discussion via the $(2l-1)$-dimensional space $\Pi_{2l-1} =\Pi_1\cap\Pi_2$. The intersection of a $(2l-1)$-dimensional space with the non-singular hyperbolic quadric $\Qplus{2l+1}$ in $\PG(2l+1,q)$ is either:   (1) a non-singular hyperbolic quadric $\Qplus{2l-1}$, (2) a cone  $L\Qplus{2l-3}$, (3)  a cone $
P\Q{2l-2}$, or (4)  a non-singular elliptic quadric $\Qmin{2l-1}$.

\begin{enumerate}
\item Let $\PG(2l-1,q)$ be an $(2l-1)$-dimensional space intersecting $\Qplus{2l+1}$ in a non-singular $(2l-1)$-dimensional 
hyperbolic quadric $\Qplus{2l-1}$. Then $\PG(2l-1,q)$ is the polar space of a bisecant line to $\Qplus{2l+1}$. Then $\PG(2l-1,q)$ lies in two tangent hyperplanes to $\Qplus{2l+1}$ and in $q-1$ hyperplanes intersecting $\Qplus{2l+1}$ in a non-singular parabolic quadric $\Q{2l}$. 
\item Let $\PG(2l-1,q)$ be an $(2l-1)$-dimensional space intersecting $\Qplus{2l+1}$ in a singular  quadric $L\Qplus{2l-3}$, then $\PG(2l-1,q)$ lies in the tangent hyperplanes to $\Qplus{2l+1}$ in the $q+1$ points $P$ of $L$.
\item
Let $\PG(2l-1,q)$ be an $(2l-1)$-dimensional space intersecting $\Qplus{2l+1}$ in a singular  quadric $P\Q{2l-2}$, then $\PG(2l-1,q)$ lies in the tangent hyperplane to $\Qplus{2l+1}$ in $P$, and in $q$ hyperplanes intersecting $\Qplus{2l+1}$ in non-singular parabolic quadrics $\Q{2l}$.
\item Let $\PG(2l-1,q)$ be an $(2l-1)$-dimensional space intersecting $\Qplus{2l+1}$ in a non-singular $(2l-1)$-dimensional elliptic quadric $\Qmin{2l-1}$, then $\PG(2l-1,q)$ lies in $q+1$ hyperplanes  intersecting $\Qplus{2l+1}$ in non-singular parabolic quadrics $\Q{2l}$.

\end{enumerate}

In the next tables, $\Qplus{2l-1}$ and $\Qmin{2l-1}$ denote non-singular hyperbolic and elliptic quadrics in $\PG(2l-1,q)$, $P\Q{2l-2}$ denotes a singular quadric with vertex the point $P$ and base a non-singular parabolic quadric in $\PG(2l-2,q)$, $L\Qplus{2l-3}$ denotes a singular quadric with vertex the line $L$ and base a non-singular hyperbolic quadric in $\PG(2l-3,q)$, $\Q{2l}$ denotes a non-singular parabolic quadric in $\PG(2l,q)$, and $P\Qplus{2l-1}$ denotes a singular quadric with vertex the point $P$ and base a non-singular hyperbolic quadric in $\PG(2l-1,q)$.

In Table 1, we denote the different possibilities for the intersection of $\Qplus{2l+1}$ with the union of two hyperplanes. We describe these possibilities by giving the formula for calculating the size of the intersection. We mention the sizes of the two quadrics which are the intersection of $\Pi_1$ and $\Pi_2$ with $\Qplus{2l+1}$, and we subtract the size of the quadric which is the intersection of $\Pi_{2l-1}=\Pi_1\cap \Pi_2$ with $\Qplus{2l+1}$.
\\

\hspace{1cm}
\begin{tabular}{|c|c|c|c|}
\hline
 & & $\Pi_{2l-1}\cap \Qplus{2l+1}$&$|\Qplus{2l+1}\cap (\Pi_1\cup\Pi_2)|$  \\
\hline
(1)& (1.1) & \Qplus{2l-1}&  $2|\Q{2l}|-|\Qplus{2l-1}|$ \\
\hline
& (1.2) & \Qplus{2l-1}& $|P\Qplus{2l-1}|+|\Q{2l}|-|\Qplus{2l-1}|$ \\ \hline
& (1.3) & \Qplus{2l-1}& $2|P\Qplus{2l-1}|-|\Qplus{2l-1}|$\\ \hline
(2) & (2.1) &$L\Qplus{2l-3}$ & $2|P\Qplus{2l-1}|-|L\Qplus{2l-3}|$ \\
\hline
(3) & (3.1)& $P\Q{2l-2}$ & $2|\Q{2l}|-|P\Q{2l-2}|$ \\
\hline
 & (3.2) &$P\Q{2l-2}$ & $|\Q{2l}|+|P\Qplus{2l-1}|-|P\Q{2l-2}|$ \\ \hline
 (4) & (4.1) & $\Qmin{2l-1}$ & $2|\Q{2l}|-|\Qmin{2l-1}|$\\ 
\hline
\multicolumn{4}{c}{   }\\
\multicolumn{4}{c}{Table 1}
\end{tabular}
\newline

We now give the sizes of these intersections of $\Qplus{2l+1}$ with the union of two hyperplanes.\\

\begin{tabular}{|c|c|c|}
\hline
 & & $|\Qplus{2l+1}\cap (\Pi_1\cup\Pi_2)|$\\
\hline
(1)& (1.1) & $2q^{2l-1}+q^{2l-2}+\cdots+q^l+q^{l-2}+\cdots+q+1$ \\ \hline
& (1.2) & $2q^{2l-1}+q^{2l-2}+\cdots+q^{l+1}+2q^l+q^{l-2}+\cdots+q+1$ \\ \hline
& (1.3) & $2q^{2l-1}+q^{2l-2}+\cdots+q^{l+1}+3q^l+q^{l-2}+\cdots+q+1$ \\ \hline
(2)& (2.1) & $2q^{2l-1}+q^{2l-2}+\cdots+q^{l+1}+2q^l+q^{l-1}+\cdots+q+1 $\\ \hline
(3)& (3.1) & $2q^{2l-1}+q^{2l-2}+\cdots+q^l+q^{l-1}+\cdots+q+1 $\\ \hline
& (3.2) & $2q^{2l-1}+q^{2l-2}+\cdots+q^{l+1}+2q^l+q^{l-1}+\cdots+q+1$ \\ \hline
& (4.1) & $2q^{2l-1}+q^{2l-2}+\cdots+q^l+2q^{l-1}+q^{l-2}+\cdots+q+1$ \\ 
\hline
\multicolumn{3}{c}{   }\\
\multicolumn{3}{c}{Table 2}
\end{tabular}
\newline

We now present in the next table the  weights of the corresponding codewords of $C_2(\Qplus{2l+1})$, and the numbers of  codewords having these weights.
\\

\begin{tabular}{|c|c|c|}
\hline
  & Weight & Number of codewords \\
  \hline
(1.3)& $w_1=q^{2l}-q^{2l-1}-q^l+q^{l-1}$ & $\frac{(q^{3l}+q^{2l})(q^{l+1}-1)}{2}$ \\ 
 \hline
(2.1)+(3.2) & $w_1+q^l-q^{l-1} $ & $\frac{(q^{2l+1}-q)(q^{l+1}-1)(q^{l-1}+1)}{2(q-1)}$+ $(q^{3l-1}-q^{l-1})(q^{l+2}-q)$\\ \hline
(1.2) & $w_1+q^{l}$ & $(q^{3l}+q^{2l})(q^{l+1}-1)(q-1)$\\
\hline
(4.1) & $w_1+2q^l-2q^{l-1} $& $\frac{q^{2l+1}(q^{l+1}-1)(q^l-1)(q-1)}{4}$ \\ \hline
 (3.1)& $w_1+2q^l-q^{l-1} $& $\frac{(q^{3l-1}-q^{l-1})(q^{l+1}-1)(q^2-q)}{2}$\\ \hline
 (1.1)& $w_1+2q^{l}$& $\frac{(q^{3l}+q^{2l})(q^{l+1}-1)(q^2-3q+2)}{4}$\\ \hline
\multicolumn{3}{c}{   }\\
\multicolumn{3}{c}{Table 3}
\end{tabular}

\begin{remark}
In the case that $q=2$, we have that the third weight coincides with the fourth. So in that special case there are only five  different weights.
\end{remark}

\subsection{The elliptic quadric in PG$(2l+1,q)$}
We know that the largest intersections of a non-singular elliptic quadric $\Qmin{2l+1}$ in $\PG(2l+1,q)$ with the other quadrics are the intersections of $\Qmin{2l+1}$ with the quadrics which are the union of two hyperplanes $\Pi_1$ and $\Pi_2$. We now discuss all the different possibilities for the intersections of $\Qmin{2l+1}$ with the union of two hyperplanes. This then gives the five or six smallest weights of the functional codes $C_2(\Qmin{2l+1})$, and the numbers of  codewords having these weights.

We again start the discussion via the $(2l-1)$-dimensional space $\Pi_{2l-1} =\Pi_1\cap\Pi_2$. The intersection of a $(2l-1)$-dimensional space with the non-singular elliptic quadric $\Qmin{2l+1}$ in $\PG(2l+1,q)$ is either:   (1) a non-singular elliptic quadric $\Qmin{2l-1}$, (2) a cone  $P\Q{2l-2}$, (3)  a cone $
L\Qmin{2l-3}$, or (4)  a non-singular hyperbolic quadric $\Qplus{2l-1}$.

\begin{enumerate}
\item Let $\PG(2l-1,q)$ be an $(2l-1)$-dimensional space intersecting $\Qmin{2l-1}$ in a non-singular $(2l-1)$-dimensional 
elliptic quadric $\Qmin{2l-1}$. Then $\PG(2l-1,q)$ is the polar space of a bisecant line to $\Qmin{2l+1}$. Then $\PG(2l-1,q)$ lies in two tangent hyperplanes to $\Qmin{2l+1}$ and in $q-1$ hyperplanes intersecting $\Qmin{2l+1}$ in a non-singular parabolic quadric $\Q{2l}$. 
\item Let $\PG(2l-1,q)$ be an $(2l-1)$-dimensional space intersecting $\Qmin{2l+1}$ in a singular  quadric $P\Q{2l-2}$, then $\PG(2l-1,q)$ lies in the tangent hyperplane to $\Qmin{2l+1}$ in the point $P$,  and in $q$ hyperplanes intersecting $\Qmin{2l+1}$ in non-singular parabolic quadrics $\Q{2l}$.
\item
Let $\PG(2l-1,q)$ be an $(2l-1)$-dimensional space intersecting $\Qmin{2l+1}$ in a singular  quadric $L\Qmin{2l-3}$, then $\PG(2l-1,q)$ lies in the tangent hyperplane to $\Qmin{2l+1}$ in the $q+1$ points $P$ of $L$.
\item Let $\PG(2l-1,q)$ be an $(2l-1)$-dimensional space intersecting $\Qmin{2l+1}$ in a non-singular $(2l-1)$-dimensional hyperbolic quadric $\Qplus{2l-1}$, then $\PG(2l-1,q)$ lies in $q+1$ hyperplanes  intersecting $\Qmin{2l+1}$ in non-singular parabolic quadrics $\Q{2l}$.

\end{enumerate}

In the next tables, $\Qplus{2l-1}$ and $\Qmin{2l-1}$ denote non-singular hyperbolic and elliptic quadrics in $\PG(2l-1,q)$, $P\Q{2l-2}$ denotes a singular quadric with vertex the point $P$ and base a non-singular parabolic quadric in $\PG(2l-2,q)$, $L\Qmin{2l-3}$ denotes a singular quadric with vertex the line $L$ and base a non-singular elliptic quadric in $\PG(2l-3,q)$, $\Q{2l}$ denotes a non-singular parabolic quadric in $\PG(2l,q)$, and $P\Qmin{2l-1}$ denotes a singular quadric with vertex the point $P$ and base a non-singular elliptic quadric in $\PG(2l-1,q)$.

In Table 4, we denote the different possibilities for the intersection of $\Qmin{2l+1}$ with the union of two hyperplanes.
\\

\hspace{1cm}
\begin{tabular}{|c|c|c|c|}
\hline
 & &$\Pi_{2l-1}\cap \Qmin{2l+1}$ & $|\Qmin{2l+1}\cap (\Pi_1\cup\Pi_2)|$  \\
\hline
(1)& (1.1) & \Qmin{2l-1}&  $2|\Q{2l}|-|\Qmin{2l-1}|$ \\
\hline
& (1.2) & \Qmin{2l-1}& $|P\Qmin{2l-1}|+|\Q{2l}|-|\Qmin{2l-1}|$ \\ \hline
& (1.3) & \Qmin{2l-1}& $2|P\Qmin{2l-1}|-|\Qmin{2l-1}|$\\ \hline
(2) & (2.1)& $P\Q{2l-2}$ & $2|\Q{2l}|-|P\Q{2l-2}|$ \\
\hline
 & (2.2) &$P\Q{2l-2}$ & $|\Q{2l}|+|P\Qmin{2l-1}|-|P\Q{2l-2}|$ \\ \hline
(3) & (3.1) &$L\Qmin{2l-3}$ & $2|P\Qmin{2l-1}|-|L\Qmin{2l-3}|$ \\
\hline
 (4) & (4.1) & $\Qplus{2l-1}$ & $2|\Q{2l}|-|\Qplus{2l-1}|$\\ 
\hline
\multicolumn{4}{c}{   }\\
\multicolumn{4}{c}{Table 4}
\end{tabular}
\newline

We now give the sizes of these intersections of $\Qmin{2l+1}$ with the union of two hyperplanes.\\

\begin{tabular}{|c|c|c|}
\hline
 & &$|\Qmin{2l+1}\cap (\Pi_1\cup\Pi_2)|$ \\
\hline
(1)& (1.1) & $2q^{2l-1}+q^{2l-2}+\cdots+q^l+2q^{l-1}+q^{l-2}+\cdots+q+1$ \\ \hline
& (1.2) & $2q^{2l-1}+q^{2l-2}+\cdots+q^{l+1}+2q^{l-1}+q^{l-2}+\cdots+q+1$ \\ \hline
& (1.3) & $2q^{2l-1}+q^{2l-2}+\cdots+q^{l+1}-q^l+2q^{l-1}+q^{l-2}+\cdots+q+1$ \\ \hline
(2)& (2.1) & $2q^{2l-1}+q^{2l-2}+\cdots+q^{l+1}+q^l+q^{l-1}+\cdots+q+1 $\\ \hline
& (2.2) & $2q^{2l-1}+q^{2l-2}+\cdots+q^{l+1}+q^{l-1}+\cdots+q+1$ \\ \hline
(3)& (3.1) & $2q^{2l-1}+q^{2l-2}+\cdots+q^{l+1}+q^{l-1}+\cdots+q+1 $\\ \hline

(4)& (4.1) & $2q^{2l-1}+q^{2l-2}+\cdots+q^{l}+q^{l-2}+\cdots+q+1$ \\ 
\hline
\multicolumn{3}{c}{   }\\
\multicolumn{3}{c}{Table 5}
\end{tabular}
\newline

We now present in the next table the  weights of the corresponding codewords of $C_2(\Qmin{2l+1})$, and the numbers of  codewords having these weights.
\\

\begin{tabular}{|c|c|c|}
\hline
  & Weight & Number of codewords \\
  \hline
(1.1)& $w_1=q^{2l}-q^{2l-1}-q^l-q^{l-1}$ & $\frac{(q^{3l+1}+q^{2l})(q^{l}-1)(q^2-3q+2)}{4}$ \\ 
 \hline
(2.1) & $w_1+q^{l-1} $ & $\frac{(q^{2l+1}+q^{l})(q^{2l}-1)(q-1)}{2}$\\ \hline
(4.1) & $w_1+2q^{l-1}$ & $\frac{q^{2l+1}(q^{l+1}+1)(q^l+1)(q-1)}{4}$ \\ \hline
(1.2) & $w_1+q^{l}$ & $(q^{3l+1}+q^{2l})(q^{l}-1)(q-1)$\\
\hline
(2.2)+(3.1)& $w_1+q^l+q^{l-1}$ & $(q^{2l}+q^{l-1})(q^{2l}-1)q+\frac{(q^{l+2}+q)(q^{2l}-1)(q^{l-1}-1)}{2(q-1)}$\\  
\hline
 (1.3)& $w_1+2q^l $& $\frac{(q^{3l+1}+q^{2l})(q^{l}-1)}{2}$\\ \hline
\multicolumn{3}{c}{   }\\
\multicolumn{3}{c}{Table 6}
\end{tabular}

\begin{remark}
In the case that $q=2$, we have that the third weight coincides with the fourth. So in that special case there are only five  different weights.
\end{remark}

\begin{theorem}\label{congruence}
Let $\mathcal{X}$ be a non-degenerate quadric (hyperbolic or elliptic) in $\pg{2l+1}$ where $l\ge 1$. All the weights $w_i$ of the code $C_{2}(\mathcal{X})$ defined on $\mathcal{X}$ are divisible by $q^{l-1}$.
\end{theorem}
\begin{proof}
Let $F$ and $f$ be two forms of degree 2 in $2l+2$ indeterminates with $l\ge 1$ and $N$ the number of common zeros of $F$ and $f$ in $\mathbb{F}_q^{2l+2}$. By the theorem of Ax-Katz \cite[p. 85]{Katz:71},  $N$ is divisible by $q^{l-1}$ since $\frac{2l+2-(2+2)}{2}=l-1$.\\ 
On the other hand,  $F$ and $f$ are homogeneous polynomials, therefore $N-1$ is divisible by $q-1$.
Let $\mathcal{X}$ and $\mathcal{Q}$ be the projective quadrics associated to $F$ and $f$, one has $\vert \mathcal{X}\cap \mathcal{Q}\vert=\frac{N-1}{q-1}$. Let $M=\frac{N-1}{q-1}$, one has 
\begin{equation}\label{Katz1}
M=\frac{kq^{l-1}-1}{q-1}=k\frac{q^{l-1}-1}{q-1}+\frac{k-1}{q-1}=k^{\prime}q^{l-1}+\pi_{l-2} 
\end{equation}
where $k$, $k^{\prime}\in \mathbb{Z}$ and $k=k^{\prime}(q-1)+1$. 
By the theorem of Ax-Katz \cite[p. 85]{Katz:71}  again, we get that the number of zeros of the polynomial $F$ in $\mathbb{F}_q^{2l+2}$ is divisible by $q^l$, so that
\begin{equation}\label{Katz2} 
\vert\mathcal{X}\vert=\frac{tq^l-1}{q-1}=t\frac{q^l-1}{q-1}+\frac{t-1}{q-1}=t^{\prime}q^l+\pi_{l-1}
\end{equation}
where $t$, $t^{\prime}\in \mathbb{Z}$ and $t=t^{\prime}(q-1)+1$.
The weight of a codeword associated to the quadric $\mathcal{X}$ is equal to:
\begin{equation}\label{poids}
w=\vert \mathcal{X}\vert -\vert\mathcal{X}\cap\mathcal{Q}\vert=\vert\mathcal{X}\vert-M
\end{equation}
Therefore, from (\ref{Katz1}), (\ref{Katz2}), and (\ref{poids}), we deduce that $w=t^{\prime}q^l-k^{\prime}q^{l-1}+q^{l-1}$.
\end{proof}

\subsection{The parabolic quadric in PG$(2l,q)$}
We know that the largest intersections of a non-singular parabolic quadric $\Q{2l}$ in $\PG(2l,q)$ with the other quadrics are the intersections of $\Q{2l}$ with the quadrics which are the union of two hyperplanes $\Pi_1$ and $\Pi_2$. We now discuss all the different possibilities for the intersections of $\Q{2l}$ with the union of two hyperplanes. This then gives the five smallest weights of the functional codes $C_2(\Q{2l})$, and the numbers of these codewords.

We proceed as follows. We start the discussion via the $(2l-2)$-dimensional space $\Pi_{2l-2}$ which is the intersection of these
two hyperplanes $\Pi_1$ and $\Pi_2$. The intersection of a $(2l-2)$-dimensional space with the non-singular parabolic quadric $\Q{2l}$ in $\PG(2l,q)$ is either:   (1) a non-singular parabolic quadric $\Q{2l-2}$, (2)  a cone $
P\Qplus{2l-3}$, (3)  a cone $P\Qmin{2l-3}$, or (4) a cone $L\Q{2l-4}$. 

For $q$ odd, we can make the discussion via the orthogonal polarity corresponding to the non-singular parabolic quadric $\Q{2l}$. For $q$ even, we need to use another approach, since then $\Q{2l}$ has a nucleus $N$. This implies that we need to make a distinction between the $(2l-2)$-dimensional spaces $\Pi_{2l-2}$ intersecting $\Q{2l}$ in a  parabolic quadric  $\Q{2l-2}$ or a quadric $L\Q{2l-4}$, containing the nucleus $N$ of $\Q{2l}$, and those not containing the nucleus $N$ of $\Q{2l}$. In \cite{JWPH:91}, these $(2l-2)$-dimensional spaces are respectively called {\em nuclear} and {\em non-nuclear}. \\

We first discuss the case $q$ odd.

\begin{enumerate}
\item Let $\PG(2l-2,q)$ be an $(2l-2)$-dimensional space intersecting $\Q{2l}$ in a non-singular $(2l-2)$-dimensional 
parabolic quadric $\Q{2l-2}$. Then $\PG(2l-2,q)$ is the polar space of a bisecant or external line to $\Q{2l}$. In the first case, $\PG(2l-2,q)$ lies in two tangent hyperplanes to $\Q{2l}$, $(q-1)/2$ hyperplanes intersecting $\Q{2l}$ in a non-singular hyperbolic quadric $\Qplus{2l-1}$, and in $(q-1)/2$ hyperplanes intersecting $\Q{2l}$ in a non-singular elliptic quadric $ \Qmin{2l-1}$. In the second case, $\PG(2l-2,q)$ lies in $(q+1)/2$ hyperplanes intersecting $\Q{2l}$ in a non-singular hyperbolic quadric $\Qplus{2l-1}$, and in $(q+1)/2$ hyperplanes intersecting $\Q{2l}$ in a non-singular elliptic quadric $ \Qmin{2l-1}$. 
\item
Let $\PG(2l-2,q)$ be an $(2l-2)$-dimensional space intersecting $\Q{2l}$ in a singular  quadric $P\Qplus{2l-3}$, then $\PG(2l-2,q)$ lies in the tangent hyperplane to $\Q{2l}$ in $P$ and in $q$ hyperplanes intersecting $\Q{2l}$ in non-singular hyperbolic quadrics $\Qplus{2l-1}$.
\item Let $\PG(2l-2,q)$ be an $(2l-2)$-dimensional space intersecting $\Q{2l}$ in a singular  quadric $P\Qmin{2l-3}$, then $\PG(2l-2,q)$ lies in the tangent hyperplane to $\Q{2l}$ in $P$, and in $q$ hyperplanes intersecting $\Q{2l}$ in non-singular elliptic quadrics $\Qmin{2l-1}$.
\item Let $\PG(2l-2,q)$ be an $(2l-2)$-dimensional space intersecting $\Q{2l}$ in a singular  quadric $L\Q{2l-4}$, then $\PG(2l-2,q)$ lies in the tangent hyperplanes to $\Q{2l}$ in the $q+1$ points $P$ of $L$.
\end{enumerate}

In the next tables, $\Q{2l-2}$ and $\Q{2l}$ denote non-singular parabolic quadrics in $\PG(2l-2,q)$ and in $\PG(2l,q)$, $\Qplus{2l-1}$ denotes a non-singular hyperbolic quadric in $\PG(2l-1,q)$, $\Qmin{2l-1}$ denotes a non-singular elliptic quadric in $\PG(2l-1,q)$, $P\Q{2l-2}$ denotes a singular quadric with vertex the point $P$ and base a non-singular parabolic quadric in $\PG(2l-2,q)$, $P\Qplus{2l-3}$ denotes a singular quadric with vertex the point $P$ and base a non-singular hyperbolic quadric in $\PG(2l-3,q)$, $P\Qmin{2l-3}$ denotes a singular quadric with vertex the point $P$ and base a non-singular elliptic quadric in $\PG(2l-3,q)$, and $L\Q{2l-4}$ denotes a singular quadric with vertex the line $L$ and base a non-singular parabolic quadric in $\PG(2l-4,q)$.

In Table 7, we denote the different possibilities for the intersection of $\Q{2l}$ with the union of two hyperplanes.
\\

\hspace{1cm}
\begin{tabular}{|c|c|c|c|}
\hline
 & & $\Pi_{2l-2}\cap \Q{2l}$&$|\Q{2l}\cap (\Pi_1\cup\Pi_2)|$ \\
\hline
(1)& (1.1) & \Q{2l-2}&  $2|\Qplus{2l-1}|-|\Q{2l-2}|$ \\
\hline
& (1.2) & \Q{2l-2}& $|\Qplus{2l-1}|+|\Qmin{2l-1}|-|\Q{2l-2}|$ \\ \hline
& (1.3) & \Q{2l-2}& $|P\Q{2l-2}|+|\Qplus{2l-1}|-|\Q{2l-2}|$\\ \hline
& (1.4) & \Q{2l-2}& $|P\Q{2l-2}|+|\Qmin{2l-1}|-|\Q{2l-2}|$\\ \hline
& (1.5) & \Q{2l-2}& $2|\Qmin{2l-1}|-|\Q{2l-2}|$\\ \hline
& (1.6) & \Q{2l-2}& $2|P\Q{2l-2}|-|\Q{2l-2}|$\\ \hline
(2) & (2.1) &$P\Qplus{2l-3}$ & $2|\Qplus{2l-1}|-|P\Qplus{2l-3}|$ \\
\hline
 & (2.2)& $P\Qplus{2l-3}$ & $|\Qplus{2l-1}|+|P\Q{2l-2}|-|P\Qplus{2l-3}|$ \\
\hline
(3) & (3.1)& $P\Qmin{2l-3}$ & $2|\Qmin{2l-1}|-|P\Qmin{2l-3}|$ \\
\hline
 & (3.2) &$P\Qmin{2l-3}$ & $|\Qmin{2l-1}|+|P\Q{2l-2}|-|P\Qmin{2l-3}|$ \\ \hline
 (4) & (4.1) & $L\Q{2l-4}$ & $2|P\Q{2l-2}|-|L\Q{2l-4}|$\\ 
\hline
\multicolumn{4}{c}{   }\\
\multicolumn{4}{c}{Table 7}
\end{tabular}
\newline

We now give the sizes of these intersections of $\Q{2l}$ with the union of two hyperplanes.\\

\begin{tabular}{|c|c|c|}
\hline
 & & $|\Q{2l}\cap (\Pi_1\cup\Pi_2)|$ \\
\hline
(1)& (1.1) & $2q^{2l-2}+q^{2l-3}+\cdots+q^l+3q^{l-1}+q^{l-2}+\cdots+q+1$ \\ \hline
& (1.2) & $2q^{2l-2}+q^{2l-3}+\cdots+q^l+q^{l-1}+q^{l-2}+\cdots+q+1$ \\ \hline
& (1.3) & $2q^{2l-2}+q^{2l-3}+\cdots+q^l+2q^{l-1}+q^{l-2}+\cdots+q+1$ \\ \hline
& (1.4) & $2q^{2l-2}+q^{2l-3}+\cdots+q^l+q^{l-2}+\cdots+q+1$ \\ \hline
& (1.5) & $2q^{2l-2}+q^{2l-3}+\cdots+q^l-q^{l-1}+q^{l-2}+\cdots+q+1$ \\ \hline
& (1.6) & $2q^{2l-2}+q^{2l-3}+\cdots+q^l+q^{l-1}+q^{l-2}+\cdots+q+1$ \\ \hline
(2)& (2.1) & $2q^{2l-2}+q^{2l-3}+\cdots+q^l+2q^{l-1}+q^{l-2}+\cdots+q+1 $\\ \hline
& (2.2) & $2q^{2l-2}+q^{2l-3}+\cdots+q^l+q^{l-1}+q^{l-2}+\cdots+q+1$ \\ \hline
(3)& (3.1) & $2q^{2l-2}+q^{2l-3}+\cdots+q^l+q^{l-2}+\cdots+q+1 $\\ \hline
& (3.2) & $2q^{2l-2}+q^{2l-3}+\cdots+q^l+q^{l-1}+q^{l-2}+\cdots+q+1$ \\ \hline
(4)& (4.1) & $2q^{2l-2}+q^{2l-3}+\cdots+q^l+q^{l-1}+q^{l-2}+\cdots+q+1$ \\ 
\hline
\multicolumn{3}{c}{   }\\
\multicolumn{3}{c}{Table 8}
\end{tabular}
\newline

\vspace{1cm}

\begin{tabular}{|c|c|c|}
\hline
  & Weight & Number of codewords \\
  \hline
(1.1)& $w_1=q^{2l-1}-q^{2l-2}-2q^{l-1}$ & $\frac{(q^{2l}-1)q^{2l-1}(q-1)(q-3)}{16}$ \\
& & +$\frac{q^{2l-1}(q^{2l}-1)(q-1)^2}{16}$  \\ \hline
(1.3)+(2.1) & $w_1+q^{l-1} $ & $\frac{(q^{2l}-1)q^{2l-1}(q-1)}{2}$+\\
& & $\frac{q^l(q^{l-1}+1)(q^{2l}-1)(q-1)}{4}$\\ \hline
(1.2) & $w_1+2q^{l-1}$ & $\frac{(q^{2l}-1)q^{2l-1}(q-1)^2}{8}$+\\
& & +$\frac{q^{2l-1}(q^{2l}-1)(q^2-1)}{8}$\\
+(1.6)+(2.2) & & $+\frac{(q^{2l}-1)q^{2l-1}}{2}+\frac{q^l(q^{l-1}+1)(q^{2l}-1)}{2}$ \\
 +(3.2)+(4.1)& & $\frac{q^l(q^{l-1}-1)(q^{2l}-1)}{2}+\frac{(q^{2l}-1)(q^{2l-2}-1)q}{2(q-1)}$\\ \hline
 (1.4)+(3.1)& $w_1+3q^{l-1}$& $\frac{(q^{2l}-1)q^{2l-1}(q-1)}{2}+\frac{q^{l}(q^{l-1}-1)(q^{2l}-1)(q-1)}{4}$\\ \hline
  (1.5) & $w_1+4q^{l-1}$ & $\frac{(q^{2l}-1)q^{2l-1}(q-1)(q-3)}{16}+\frac{q^{2l-1}(q^{2l}-1)(q-1)^2}{16}$\\
 \hline
\multicolumn{3}{c}{   }\\
\multicolumn{3}{c}{Table 9: Weights and number of codewords for $q$ odd}
\end{tabular}
\vspace{0.5cm}

We now discuss the case $q$ even. Here Q$(2l,q)$ has a nucleus $N$.

\begin{enumerate}
\item Let $\PG(2l-2,q)$ be an $(2l-2)$-dimensional space intersecting $\Q{2l}$ in a non-singular $(2l-2)$-dimensional 
parabolic quadric $\Q{2l-2}$. If $\PG(2l-2,q)$ is non-nuclear, then $\PG(2l-2,q)$ lies in one tangent hyperplane, the hyperplane $\langle\PG(2l-2,q),N\rangle$, in $q/2$  hyperplanes intersecting $\Q{2l}$ in a non-singular hyperbolic quadric $\Qplus{2l-1}$, and in $q/2$ hyperplanes intersecting $\Q{2l}$ in a non-singular elliptic quadric $ \Qmin{2l-1}$. If $\PG(2l-2,q)$ is nuclear, then $\PG(2l-2,q)$ lies in $q+1$ tangent hyperplanes to $\Q{2l}$.
\item
Let $\PG(2l-2,q)$ be an $(2l-2)$-dimensional space intersecting $\Q{2l}$ in a singular  quadric $P\Qplus{2l-3}$, then $\PG(2l-2,q)$ lies in the tangent hyperplane to $\Q{2l}$ in $P$, and in $q$ hyperplanes intersecting $\Q{2l}$ in non-singular hyperbolic quadrics $\Qplus{2l-1}$.
\item Let $\PG(2l-2,q)$ be an $(2l-2)$-dimensional space intersecting $\Q{2l}$ in a singular  quadric $P\Qmin{2l-3}$, then $\PG(2l-2,q)$ lies in the tangent hyperplane to $\Q{2l}$ in $P$, and in $q$ hyperplanes intersecting $\Q{2l}$ in non-singular elliptic quadrics $\Qmin{2l-1}$.
\item Let $\PG(2l-2,q)$ be an $(2l-2)$-dimensional space intersecting $\Q{2l}$ in a singular  quadric $L\Q{2l-4}$, then $\PG(2l-2,q)$ lies in the tangent hyperplanes to $\Q{2l}$ in the $q+1$ points $P$ of $L$.
\end{enumerate}

In Table 7, we denoted the different possibilities for the intersection of $\Q{2l}$ with the union of two hyperplanes, and in Table 8, the corresponding sizes for the intersections. We now present in Table 10 the number of codewords having the corresponding weights.
\newline

\begin{tabular}{|c|c|c|}
\hline
  & Weight & Number of codewords \\
  \hline
(1.1) & $w_1=q^{2l-1}-q^{2l-2}-2q^{l-1}$ & $\frac{(q^{2l}-1)q^{2l-1}(q-2)(q-1)}{8}$ \\ \hline
(1.3)+(2.1)  & $w_1+q^{l-1} $ & $\frac{(q^{2l}-1)q^{2l-1}(q-1)}{2}$+ $\frac{q^l(q^{l-1}+1)(q^{2l}-1)(q-1)}{4}$\\ \hline
(1.2)+(1.6) & $w_1+2q^{l-1}$ & $\frac{(q^{2l}-1)q^{2l}(q-1)}{4}$+ $\frac{q^{2l-1}(q^{2l}-1)}{2}+$\\
+(4.1) & & $\frac{q(q^{2l-2}-1)(q^{2l}-1)}{2(q-1)}+$\\
+(2.2)+(3.2)& & $\frac{q^{l}(q^{l-1}+1)(q^{2l}-1)}{2}+\frac{q^{l}(q^{l-1}-1)(q^{2l}-1)}{2}$ \\
  \hline
 (1.4)+(3.1)& $w_1+3q^{l-1}$& $\frac{(q^{2l}-1)q^{2l-1}(q-1)}{2}+\frac{q^{l}(q^{l-1}-1)(q^{2l}-1)(q-1)}{4}$\\ \hline
  (1.5) & $w_1+4q^{l-1}$ & $\frac{(q^{2l}-1)q^{2l-1}(q-1)(q-2)}{8}$\\
 \hline
\multicolumn{3}{c}{   }\\
\multicolumn{3}{c}{Table 10: Weights and number of codewords for $q$ even}
\end{tabular}

\begin{theorem}
Let $\mathcal{X}$ be a non-degenerate parabolic quadric in $\pg{2l}$ where $l\ge 1$. All the weights $w_i$ of the code $C_{2}(\mathcal{X})$ defined on $\mathcal{X}$ are divisible by $q^{l-1}$.
\end{theorem}
\textbf{Proof:}
It is analogous to the one of Theorem \ref{congruence}.


Address of the authors:\\
 CNRS, Institut de Math\'ematiques de Luminy, Luminy Case 907, 13288 Marseille Cedex 9, France.

\noindent
F. A. B. Edoukou: edoukou@iml.univ-mrs.fr,\\
F. Rodier: rodier@iml.univ-mrs.fr, http://iml.univ-mrs.fr/$\sim$rodier/\\

 Department of pure mathematics and computer algebra,           			
                         Ghent University,      
                         Krijgslaan 281-S22,
                        9000 Ghent,
                        Belgium.
                        
 \noindent 
A. Hallez: athallez@cage.ugent.be,\\
L. Storme: ls@cage.ugent.be, http://cage.ugent.be/$\sim$ls
\end{document}